\documentstyle[twoside]{article}
\begin{document}
\pagestyle{myheadings}
\markboth{E. Vassiliou}{ASSOCIATED VECTOR SHEAVES}

\newtheorem{thm}{Theorem}
\newtheorem{lem}{Lemma}
\newtheorem{cor}{Corollary}
\newtheorem{Def}{Definition}
\newtheorem{rem}{Remark}

\newcommand{\lra}{\longrightarrow}
\newcommand{\prt}{\partial}
\newcommand{\sst}{\subseteq}
\newcommand{\bm}[1]{\mbox{\boldmath$#1$}}

\newcommand{\ca}{{\cal A}}
\newcommand{\cb}{{\cal B}}
\newcommand{\cC}{{\cal C}}
\newcommand{\ce}{{\cal E}}
\newcommand{\cg}{{\cal G}}
\newcommand{\cl}{{\cal L}}
\newcommand{\cp}{{\cal P}}
\newcommand{\cu}{{\cal U}}

\newcommand{\ua}{U_{\alpha}}
\newcommand{\ub}{U_{\beta}}
\newcommand{\uab}{U_{\alpha \beta}}
\newcommand{\uua}{U \cap \ua}
\newcommand{\uub}{U \cap \ub}
\newcommand{\uuab}{U \cap \ua \cap \ub}

\newcommand{\ad}{{\ca d}}
\newcommand{\ga}{g_{\alpha}}
\newcommand{\gb}{g_{\beta}}
\newcommand{\gab}{g_{\alpha \beta}}

\newcommand{\ha}{h_{\alpha}}
\newcommand{\hb}{h_{\beta}}
\newcommand{\fa}{\phi_{\alpha}}
\newcommand{\fb}{\phi_{\beta}}
\newcommand{\sa}{s_{\alpha}}
\newcommand{\sbe}{s_{\beta}}
\newcommand{\Sa}{\Sigma_{\alpha}}
\newcommand{\Sbe}{\Sigma_{\beta}}
\newcommand{\Pa}{\Pi_{\alpha}}
\newcommand{\Pb}{\Pi_{\beta}}
\newcommand{\gr}{\varrho}
\newcommand{\gra}{\alpha}
\newcommand{\grb}{\beta}
\newcommand{\gre}{\varepsilon}
\newcommand{\vth}{\varTheta}
\newcommand{\vf}{\varphi}
\newcommand{\of}{\overline{\vf}}

\newcommand{\wa}{\omega_{\alpha}}
\newcommand{\wb}{\omega_{\beta}}
\newcommand{\adg}{\ad(\gab^{-1})}
\newcommand{\rg}{\rho(\gab^{-1})}
\newcommand{\cau}{\ca(U)}
\newcommand{\caua}{\ca(\ua)}

\newcommand{\ma}{M_{n}(\ca)}
\newcommand{\mau}{M_{n}(\cau)}
\newcommand{\maua}{M_{n}(\caua)}

\newcommand{\gla}{{\cal GL}(n,\ca)}

\newcommand{\glau}{{\cal GL}(n,\ca)(U)}

\newcommand{\W}{\Omega^1}
\newcommand{\WL}{\Omega^1\otimes_{\ca}\cl}
\newcommand{\Wa}{\Omega^{1}\otimes_{\ca}\ca}
\newcommand{\Wma}{\Omega^{1}\otimes_{\ca}\ma}
\newcommand{\Wmau}{\Omega^{1}(U)\otimes_{\cau}\mau}
\newcommand{\Wmaua}{\Omega^{1}(\ua)\otimes_{\caua}\maua}
\newcommand{\WWma}{\Omega^{2}\otimes_{\ca}\ma}
\newcommand{\We}{\W\otimes_{\ca}\ce}

\title{VECTOR SHEAVES ASSOCIATED WITH PRINCIPAL SHEAVES $^*$}
\author{Efstathios Vassiliou}
\date{}
\maketitle
\thispagestyle{empty}

\begin{abstract}
\noindent
In the framework of Abstract Differential Geometry, especially
that dealing with vector sheaves (as expounded in \cite{Mal3}) and
principal sheaves (initiated by \cite{V1}), we show that to a given
principal sheaf $(\cp,\cg,X,\pi)$ together with a representation
$\vf : \cg \rightarrow \gla$, we associate a vector sheaf $(\ce,X,p)$.
If $\vf$ is compatible with the representations of $\cg$ and $\gla$
into appropriate sheaves of Lie algebras, as well as with  the
Maurer-Cartan (or logarithmic) differentials of  the same sheaves
of groups, then every connection on $\cp$ induces an $\ca$-linear
connection on $\ce$. An example is provided by the principal sheaf of
frames of a vector sheaf.
\end{abstract}



\medskip
\noindent
{\bf *} This paper is in final form and no version of it will be
           submitted for publication elsewhere

\setcounter{section}{-1}
\section{Introduction}   
Gauge theories are, roughly speaking, built on principal bundles
$(P,G,X,\pi)$ and their connections. This is a consequence of the
fact that observations and measurements in physics lead to certain
sections of a parametrized group, in general non abelian. However,
Lie groups and principal bundles are quite complicated objects and
one is looking for a reduction of the non-commutative framework to
a commutative one, the latter being described by a vector bundle.
This can be often achieved by an appropriate representation of 
$G$ into a vector space (in this respect we refer to \cite{Bl}).

The aim of this note is to examine the analogous situation in the
context of {\em Abstract Differential Geometry\/}. As a matter of
fact, the present author has initiated a research program devoted 
to the geometry of principal sheaves (see \cite{V1}--\cite{V4})
influenced by the geometry of vector sheaves expounded in \cite{Mal3}.
These abstractions are developed in a completely
{\em algebrotopological\/} setting, without any differentiability,
in spite of the wide use of the adjective ``differential''
accompanying various terms in order to remind the analogy with the
classical geometry of ordinary (smooth) fiber bundles.

In the present abstract approach we consider a principal sheaf
$(\cp,\cg,X,\pi)$ and a representation of the form
$\vf : \cg \lra \gla$, where $\ca$ is a sheaf of unital,
commutative and associative algebras. Thus $\ca^n$ is the vector sheaf
in which the structure sheaf of groups $\cg$ is represented.
We show that such a representation leads to a vector sheaf
$(\ce,X,p)$ associated with $\cp$ (Section 2). In the
sequel (Section 3), under some additional assumptions pertaining
to the compatibility of $\vf$ with the Maurer-Cartan (or logarithmic)
differentials of $\cp$ and $\gla$, as well as with the representations
of the latter into certain sheaves of Lie algebras, we prove that the
connections on $\cp$ (in the sense of \cite{V1}) induce
$\ca$-connections on $\ce$ (in the sense of \cite{Mal3}).
The converse is not always true unless extra conditions are imposed on
$\vf$. An example is provided by the principal {\em sheaf of frames\/}
of a given vector sheaf (already studied in \cite{V2}), in which case 
we have the trivial representation of $\gla$.

Since the notations and terminology used throughout are not yet
standard, the preliminary Section 1 contains a brief account of
the material essentially needed,in order to make the note as
self sufficient as possible, referring for details to the relevant
literature.

\section{Preliminaries}  
{\bf 1.} Our setting is based on a fixed {\em algebraized space\/}
{ $(X,\ca)$, where $X$ is a topological space and $\ca$ a sheaf
(over $X$) of {\em unital, commutative\/} and {\em associative\/}
$\bm{K}$-{\em algebras\/} ($\bm{K}=\bm{R},\bm{C}$). For instance,
in the classical case of a real smooth manifold $X$, we take
$\ca=\cC_X^{\infty}$, the sheaf of germs of smooth functions on X.
For other examples we refer to \cite[Chapter 10]{Mal3}.

To such an algebraized space we also attach a {\em differential
triad\/} $(\ca, d, \W)$, where $\W$ is an $\ca$-module (over $X$) and 
$d:\ca \lra\W$ a derivation of $\ca$; that is, a
$\bm{K}$-linear morphism  satisfying the {\em Leibniz condition\/}
\[
        d(s\cdot t)=s\cdot d(t)+t\cdot d(s),
\]
for any (local) sections $s,t \in \ca(U)$ and $U\sst X$ open.
Note that in the previous formula we have identified a sheaf with
the sheaf of germs of its sections, a convenient fact which will be
often used below.

In the classical case, $\W$ is nothing but the sheaf of germs of smooth
1-forms on $X$. In the abstract (algebrotopological) framework we are
dealing with, differential triads always exist by K\"{a}hler's theory
of differentials (for details \cite{Mal2}, \cite[Chapter 11,
Sections 5--6]{Mal3}).

\medskip
\noindent
{\bf 2.} Among the objects of prime interest here are principal
sheaves, originally considered (in a different context) by
A.~Grothendieck \cite{Gr}. More precisely, a {\em principal sheaf\/}
over $X$ is described by a quadruple $\cp \equiv (\cp,\cg,X,\pi)$,
where $\pi$ is the projection of $\cp$ on $X$ and $\cg$ is a sheaf
of groups representing simultaneously the {\em structure sheaf\/}
and the {\em structural type\/} of $\cp$. This means that there exists
a (right) action $\cp \times _X \cg \lra \cp$, as well as an open
covering $\cu = \{\ua|\gra \in I\}$ of $X$ together with local
$\cg$-equivariant isomorphisms $\phi_{\gra} : \cp|_{\ua}
\lra \cg|_{\ua}$.

However, in order to built up an abstract differential geometry on
$\cp$, in particular a gauge theory, we enrich the structure sheaf
with two additional properties. In fact, we assume that $\cg$ is a
sheaf of groups of {\em Lie-type\/}, by which we mean that:

i) There exists a {\em representation\/} (: a continuous morphism
of sheaves of groups) $\gr: \cl \lra \ca ut(\cl)$ of $\cg$ in an
$\ca$-module of Lie algebras $\cl$;

ii) There exists a morphism (of sheaves of sets)
          $\prt : \cg \lra \WL$,
called {\em Maurer-Cartan\/ {\em or} logarithmic differential\/},
such that
\[
\prt(s \cdot t) = \gr(t^{-1}) .\prt(s) + \prt(t),
\]
for every $s,t \in \cg(U)$ and $U \sst X$ open.
The first term of the right-hand side of the previous formula denotes
the result of the natural action of $\cg$ on $\WL$ induced by $\gr$.
To be more explicit, for any $g \in \cg$ and any decomposable
element $\omega \equiv \theta \otimes u \in \WL$, we set
\begin{equation}    
  \varrho(g).\omega \equiv (1\otimes \varrho(g)).\omega
         := \theta \otimes \varrho(g)(u),
\end{equation}
where $1$ here denotes the identity of $\Omega^1$. We extend this
action by linearity to arbitrary elements.

$\cp$ admits a family of {\em natural (local) sections\/}
\[
  \sa := \psi_{\gra} \circ \mbox{\bf 1}|_{\ua} \in \cp;
                           \quad \gra \in I,
\]
where {\bf 1} is the unit section of $\cg$ (: {\bf 1}(x) is the
unit of the fiber $\cg_x$). 

As an example we take the sheaf $\cp$ of germs of smooth sections
of a principal fibre bundle $(P,G,X,p)$. It is a principal sheaf with
structure sheaf $\cg$ the sheaf of germs of smooth $G$-valued maps on
$X$. $\cg$ is of Lie-type with $\cl$ being now the sheaf of germs of
smooth maps on $X$ with values in the Lie algebra of $G$. In this case
$\gr$ and $\partial$ are obtained by the sheafification of the adjoint
representation and the total (logarithmic) differential respectively.
For complete details we refer e.g to \cite{V1, V4}

\medskip
\noindent
{\bf 3.} A typical abstract example of a sheaf of groups of Lie-type,
which will play an important role in the sequel, is the sheaf
$\gla$ generated by the complete presheaf of groups
$U \mapsto GL(n,\ca(U))$, U running in the topology of $X$. Hence,
\begin{equation}      
     \gla(U) \cong GL(n,\ca(U)) \cong Lis_{\ca|_U}(\ca^n|_U,\ca^n|_U).
\end{equation}    
Now $\cl \equiv \ma$, the sheaf generated by the complete presheaf
of Lie algebras $U \mapsto \mau$; thus
\begin{equation}   
  \ma(U) \cong \mau \cong \ca^{n^2}(U),
\end{equation}
for every open $U \sst X$.

There exists an (adjoint) representation
            $\ca d : \gla \lra \ca ut(\ma)$
obtained as follows: Let $U$ be any open subset of $X$. We define
the morphism of sections
\[
 \ca d^U : \gla(U) \lra \ca ut((\ma))(U) \cong Aut(\ma|_U, \ma|_U)
\]
by requiring that, for any $g \in \glau$, $\ca d^U(g)$ to be the
automorphism generated by the automorphisms
of presheaves
\[
   \big (\ca d^U(g)\big)_V : \ma(V) \lra \ma(V):\  a \mapsto
        g\cdot a \cdot g^{-1};\quad a \in \ma(V),
\]
for all open $V \sst U$, with the identifications (2) and (3) being
applied here.

The corresponding Maurer-Cartan differential
$\widetilde{\partial} : \cg \lra \Wma$ is given by
   $\partial(a) := a^{-1}\cdot\bm{d}(a)$,
for every $a \in \mau$ and $U \sst X$ open, where
$\bm{d} : \gla \lra \Wma$ is the extension of $d$ (of the initial
differential triad); i.e.,
  $\bm{d}(a) := (da_{ij})$, for every $a = (a_{ij}) \in GL(n,\ca(U))$. 

\medskip
\noindent
{\bf 4.} The last fundamental notion immediately needed in the next
section is that of a {\em vector sheaf\/}. This is a sheaf
$\ce \equiv (\ce,X,p)$ which is a locally free $\ca$-module (over X).
Hence, there exist an open cover, say, ${\cal U} = \{U_{\gra}|\gra \in
I\}$ and $\ca|_{\ua}$-isomorphisms $\psi_{\gra} : \ca^n|_{\ua} \lra
\ce|_{\ua}$. The complete study of vector sheaves and their geometry
is the content of \cite{Mal3}.

\section{Associated sheaves}   
In this section we fix a principal sheaf $\cp \equiv (\cp,\cg,X,\pi)$
and a representation of the form
 \[
     \vf : \cg \lra \gla.
 \]
We shall construct a vector sheaf of rank n, associated with $\cp$.
To this end, for each open $U \sst X$, we consider the quotient set
$Q(U) := \cp(U)\times\ca^n(U)\,/\,\cg(U)$ determined by the
equivalence relation
\[
  (s,a) \sim (t,b) \quad \Longleftrightarrow \quad \exists!\;
  g\in\cg(U)\,: \, t = s\cdot g\, , \, b = \vf(g^{-1})\cdot a,
\]
for every $s,t \in \cp(U)$ and $a,b \in \ca^n(U)$.

Running now $U$ in the topology of $X$, we obtain a (not necessarily
complete) presheaf $U \mapsto Q(U)$ generating the quotient sheaf
\[
           \ce := \cp \times _X \ca^n \, / \, \cg,
\]
with base X and a projection $p$ defined in the obvious way.
This is, by definition, {\em the sheaf associated with\/} $\cp$
{\em by the representation\/} $\vf$.

With regard to the previous construction one may consult \cite{Go}.
We note that the last quotient can be also constructed, in an
equivalent way, by defining (fibrewise) on $\cp \times _X \ca^n$
an analogous (global) equivalence relation (see \cite{Gr}).

\begin{lem}      
$(\ce,X,p)$ is a sheaf locally isomorphic to $\ca^n$ with corresponding
cocycle $(G_{\gra\grb})=(\vf(g_{\gra\grb})) \in
                 Z^1(\cu,\gla)$,
where $(g_{\gra\grb}) \in Z^1(\cu,\cg)$ is the cocycle of the
principal sheaf $\cp$.
\end{lem}
\noindent
{\bf Proof.} Fix a $\ua \in \cu$. For every open $V \sst \ua$
we define the map
\[
    \psi^{\gra}_V : \ca^n(V) \ni f \mapsto [\sa|_V,f] \in Q(V),
\]
where $\sa \in \cp(\ua)$ is the natural section over $\ua$.

It is immediate that $\psi^{\gra}_V$ is 1-1. On the other hand, for
a given $[\sigma,h] \in Q(V)$, the section $\vf(g)\cdot h \in \ca^n(V)$,
with $g$ determined by $\sigma = \sa|_V \cdot g$, gives that
$\psi^{\gra}_V (\vf(g)\cdot h) = [\sigma,h]$, which implies that
 $\psi^{\gra}_V$ is onto. In this way we obtain a morphism
 $\{\psi^{\gra}_V\}_{V\sst\ua}$ between the presheaves
 $V \mapsto \ca^n(V)$ and $V \mapsto Q(V)$, generating an
isomorphism (of sheaves of sets)
$\psi_{\gra} : \ca^n|_{\ua} \stackrel{\simeq}{\lra} \ce|_{\ua}$.
This shows the first claim of the statement.

By definition, $G_{\gra\grb}= \psi_{\gra}^{-1} \circ \psi_{\grb}$,
where now both the isomorphisms are restricted on appropriate
sheaves over $\uab:= \ua \cap \ub$ (for simplicity we omit explicit
expressions like $\psi_{\gra}|_{\uab}$). Hence, $G_{\gra\grb}$
is generated by $(\psi_V^{\gra})^{-1} \circ \psi_V^{\grb}$, for all
open $V \sst \uab$. As a result, for every $h \in \ca^n(V)$, we
check that
\[
  \big((\psi_V^{\gra}\big)^{-1} \circ \psi_V^{\grb})(h) =
      (\psi_V^{\gra})\big([\sbe|_V,h]\big) = \vf(\gab|_V)\cdot h.
\]
Using the identification (2), we obtain
 $(\psi_V^{\gra})^{-1} \circ \psi_V^{\grb} = \vf(\gab|_V)$.
We prove the second claim by taking all open $V \sst \uab$.
\hfill $\Box$
  
\begin{thm}    
$\ce \equiv (\ce,X,p)$ is a vector sheaf (of rank n).
\end{thm}

\noindent
{\bf Proof.} Each isomorphism (of sheaves of sets)  
 $\psi_{\gra} : \ca^n|_{\ua} \lra \ce|_{\ua}$ induces (fibrewise)
 on $\ce|_{\ua}$
  the operations
 \[
   \Sa : \ce|_{\ua}\times_{\ua} \ce|_{\ua} \lra \ce|_{\ua}\,;\,                 
   \Pa : \ca^n|_{\ua} \times_{\ua}\ce|_{\ua}\lra \ce|_{\ua},
 \]
respectively given by
\begin{center}
 $\Sa(u,v)\equiv u+v := \psi_{\gra}\big(\psi_{\gra}^{-1}(u)+                                   \psi_{\gra}^{-1}(v))\, ,\,
                           \psi_{\gra}^{-1}(v)\big)$
\end{center}
\begin{center}
  $\Pa(a\cdot u) \equiv a\cdot u := \psi_{\gra}\big(a\cdot
                              \psi_{\gra}^{-1}(u)\big)$, 
\end{center}
for every $u,v \in \ce_x$ and $a \in \ca_x$ with $x \in \ua$.

Since
 $\Sa=\psi_{\gra}\circ\overline{\Sigma}\circ(\psi_{\gra},
                             \psi_{\gra})$
and
 $\Pa=\psi_{\gra}\circ\overline{\Pi}\circ(\psi_{\gra},\psi_{\gra})$,
where $\overline{\Sigma}$ and $\overline{\Pi}$ are the respective
(continuous) operations of the $\ca$-module $\ca^n$,
appropriately restricted over $\ua$, it follows that $\Sa$ and $\Pa$
are also continuous morphisms giving on $\ce|_{\ua}$ the structure of
an $\ca|_{\ua}$-module such that $\psi_{\gra}$ is an
$\ca|_{\ua}$-linear isomorphism. This determines the desired local
structure of $\ce$.

The previous local operations globalize to corresponding
continuous operations on $\ce$ since $\Sa=\Sbe$ and $\Pa=\Pb$ on
the overlappings. Indeed, for any
$(u,v)\in \ce|_{\uab}\times_{\uab}\ce|_{\uab}$, using the
identification (2) and the previous Lemma, we have that
\begin{eqnarray*}
  \Sbe(u,v) & = &(\psi_{\gra}\circ G_{\gra\grb})
      \big(\psi_{\grb}^{-1}(u)+ \psi_{\grb}^{-1}(v)\big)\\
            & = & (\psi_{\gra}\circ \psi_{\grb}^{-1})(u) +
       (\psi_{\gra}\circ \psi_{\grb}^{-1})(u)\\
   & = &\psi_{\gra}\big(\psi_{\gra}^{-1}(u)+\psi_{\gra}^{-1}(v)\big) =
                          \Sa(u,v)
\end{eqnarray*}
and similarly for the multiplications. Therefore, $\ce$ becomes an
$\ca$-module.  \hfill $\Box$

\medskip
For the sake of completeness, we examine the relationship between
the (global) sections of $\ce$ and certain morphisms corresponding
to the classical tensorial maps. In fact, a morphism
(of sheaves of sets) $f : \cp \lra \ca^n$ is said to be
{\em tensorial\/} if
 \[
   f(s\cdot g) = \vf(g^{-1}) \cdot f(s)\,; \quad (s,g) \in
                            \cp(U)\times \cg(U),
 \]
for every open $U \sst X$. Clearly, the product of the right-hand side
is well defined by the obvious action of $\gla$ on the left of $\ca^n$. 
As a result, we prove

\begin{thm}   
Tensorial morphisms $f : \cp \lra \ca^n$ correspond bijectively to 
global sections of $\cp$.
\end{thm}
\noindent
{\bf Proof.} Let $f$ be a given tensorial morphism. For a
$\ua \in {\cal U}$, we set $\sigma_{\gra} := [\sa,f(\sa)]$
(recall that $\sa$ is the natural section of $\cp$ over $\ua$ and $f$
is now the induced morphism of sections). Since $\sigma_{\gra} \in
(\cp(\ua) \times \ca^n(\ua)/ \sim) \subset \ce(\ua)$, we obtain
a family of local sections $(\sigma_{\gra})$ of $\ce$. However, over
$\uab$, we have that
\[
  \sigma_{\grb} = [\sa \cdot \gab,\vf(\gab^{-1})\cdot
   f(\sa)] = [\sa,f(\sa)] = \sigma_{\gra};
\]
hence we can define a global section $\sigma \in \ce(X)$ by setting
$\sigma|_{\ua} := \sigma_{\gra}$.                     

Conversely, let $\sigma \in \ce(X)$ be given a section. For an open
$U\sst X$, we define the map $f_U : \cp(U) \lra \ca^n(U)$ by
requiring that
\begin{equation}  
  f_U(s)|_{\uua} := \vf(\ga^{-1})\cdot\psi_{\gra}^{-1}
                                       (\sigma|_{\uua}),
\end{equation}
for every $s \in \cp(U)$ and with $\ga \in \cg(\ua)$ determined by
$s|_{\uua} = \sa|_{\uua} \cdot \ga$. We check that $f_U$ is defined
by gluing the restrictions given by (4), for all $\ua \in \cu$. Indeed,
for $\ub \in \cu$, we have the analogous expression
\begin{equation}   
  f_U(s)|_{\uub} := \vf(\gb^{-1})\cdot\psi_{\grb}^{-1}(\sigma|_{\uub}),    
\end{equation}
with $\gb \in \cg(\ub)$ satisfying $s|_{\uub} = \sbe|_{\uub}
                                                   \cdot \gb$.
Therefore, over $\uuab$, $\ga = \gab \cdot \gb$ . Omitting, for
simplicity the explicit mention of the restrictions on $\uuab$ of
the sections involved, we see that (see also Lemma 1)
\[
 \vf(\gb^{-1})\cdot\psi_{\grb}^{-1}(\sigma) =
  \vf(\ga^{-1}\cdot \gab)\cdot G_{\beta\alpha}\cdot
                                    \psi_{\gra}^{-1}(\sigma)=
         \vf(\ga^{-1}) \cdot \psi_{\gra}^{-1}(\sigma),
\]
which proves that (4) and (5) coincide on $\uuab$ and $f_U$ is well
defined by the gluing process.

Finally, for any $s \in \cp(U)$ and $g \in \cg(U)$, we have that
\[
  f_U(s\cdot g)|_{\uua} = \vf(g|_{\uua} \cdot \ga^{-1})\cdot
             \psi_{\gra}^{-1}(\sigma|_{\uua})=
            \vf(g^{-1})|_{\uua} \cdot f_U(s)|_{\uua},
\]
for every $\ua \in \cu$; thus $f_U(s\cdot g) = \vf(g^{-1})
 \cdot f_U(s)$. Varying $U$ in the topology of $X$, we obtain a
morphism of presheaves generating a tensorial morphism $f$ and the
proof is now complete. \hfill   $\Box$

\begin{rem}
{\em In all the previous construction it is not necessary to
assume that $\cg$ is a sheaf of groups of Lie-type (see Paragraph
1.2), a fact which will be needed in the study of connections below.}
\end{rem}

\section{Connections on associated sheaves}  
In this section we consider a principal sheaf $\cp$ with structure
sheaf $\cg$ of Lie-type. We recall that (see \cite{V1}) a {\em
connection\/} on $\cp$ (or gauge potential, in the terminology of
\cite{Bl}) is a morphism of sheaves of sets $D: \cp \lra \WL$
satisfying
\begin{equation}                    
 D(s\cdot g) = \rg.D(s) + \prt(g) \,,
\end{equation}
for any $s \in \cp (U)$, $g \in \cg(U)$ and $U \subseteq X$ open.

A connection $D$ is equivalently determined by the family
of local sections
\[
  \wa := D(\sa) \in (\WL)(\ua); \quad \gra \in I,
\]
which are called, following the classical terminology, {\em the local
connection forms\/} (or local gauge potentials) of $D$. They
satisfy the (local) gauge transform
\begin{equation}                 
  \wb = \rg.\wa + \partial(\gab)
\end{equation}
on each $\uab \ne \emptyset$ (see \cite[Theorem 5.4]{V1}).

On the other hand (see \cite[Vol. II, Chapter 6, Section 3]{Mal3}),
an {\em $\ca$-connection\/} on a vector sheaf $\ce$ (of rank $n$)
is a $\bm K$-linear morphism $ \nabla : \ce \lra \ce \otimes_{\ca}\W$
satisfying the {\em Leibniz-Koszul\/} condition
\begin{equation}       
 \nabla (a \cdot s) = a \cdot \nabla (s) + s\otimes d(a),
\end{equation}
for every $a \in \ca(U)\,, s \in \ce(U)$ and $U \subseteq X$ open.

Equivalently (see also \cite[Chapter 7]{Mal3}), $\nabla$ is fully
determined by corresponding local connection forms as follows:
For each $\ua$, the $\ca(\ua)$-module $\ce(\ua)$ admits a natural
basis $e^{\gra} := (e^{\gra}_1, \ldots, e^{\gra}_n)$ with
 \[
e^{\gra}_i(x):= \psi_{\gra}(0_x, \ldots, 1_x, \ldots , 0_x);
  \qquad x \in \ua,
\]
where $0_x$ and $1_x$ (in the i-th entry) are the zero and unit
element of the stalk $\ca_x$ respectively. Evaluating now $\nabla$
on the sections of the basis, we obtain the expressions
\[
      \nabla (e^{\gra}_j) = \sum_{i=1}^n e^{\gra}_i \otimes
           \theta^{\gra}_{ij}\,; \qquad  1 \leq \;j\; \geq n,
\]
with $\theta^{\gra}_{ij} \in \Omega^1(\ua)$, forming thus a matrix
$(\theta^{\gra}_{ij}) \in M_n(\Omega^1(\ua))$, for every
$\gra \in I$. In virtue of (3), we check that
\begin{equation}  
        (\Wma)(\ua) \cong \Wmaua \cong  M_n(\Omega^1(\ua));   
\end{equation}
hence, $(\theta^{\gra}_{ij})$ can be identified with
a section $\theta_{\gra} \in (\Wma)(\ua)$. The sections
$(\theta_{\gra})_{\gra \in I}$, are the {\em local connection forms
of\/} $\nabla$ and satisfy the analog of (7), namely
\begin{equation}      
    \theta_{\grb} = \ca d(G_{\gra\grb}^{-1}).\theta_{\gra} +
             \widetilde{\partial}(G_{\gra\grb}),
\end{equation}
where $(G_{\gra\grb})$ is the cocycle of $\ce$.
This is a consequence of (8) and routine, though tedious,
calculations.

\medskip
We come now to the following basic

\begin{Def}   
{\em A representation $\vf : \cg \lra \gla$ is said to be of
{\em Lie-type\/} if there exists a morphism of sheaves of Lie
algebras $\of : \cl \lra\ma$ such that the following
conditions hold:
\begin{center}
 $ \widetilde{\partial}\circ\vf = (1\otimes \of)\circ\partial$
\end{center}
\begin{center}
  $\of\circ\gr(g) = \ca d(\vf(g)) \circ \of;
                  \quad g\in\cg,$
\end{center}
where, for simplicity, we have set $1 = id|\W$.}
\end{Def}
Clearly, the previous conditions express the compatibility of
$\vf$ and $\of$ with the Maurer-Cartan differentials of $\cg$ and
$\gla$, as well as with the their representations $\gr$ and
$\ca d$. For a more general situation see also
\cite[Definition 3.6]{V3}. Note that in the classical case $\of$
is the morphism of Lie algebras induced by the differential of
$\vf$ and the above conditions
are always true.
 
\begin{thm}   
Let $\vf : \cg \lra \gla$ be a representation of Lie-type. Then,
every connection on $\cp$ induces an $\ca$-linear connection
on the associated vector sheaf $\ce$.
\end{thm}

\noindent
{\bf Proof.} For a given connection $D \equiv (\wa)$ on $\cp$,
we set
\begin{equation}    
 \theta_{\gra} := (1\otimes \of)(\wa), \quad \gra \in I.
\end{equation}
Then, in virtue of (1), Lemma 1 and Definition 1, equality (6)
implies that
\begin{eqnarray*}
   \theta_{\grb} & = & (1\otimes\of)\big((1\otimes\gr(\gab)).\wa
          + \partial(g_{\gra\grb})\big)\\ 
 & = & \big(1\otimes\ca d(\vf(g_{\gra\grb}^{-1}))\circ \of\big).\wa
       + (\overline{\partial} \circ \vf)(g_{\gra\grb})\\ 
  & = & \big(1\otimes\ca d(G_{\gra\grb}^{-1})\big).
         \theta_{\gra}+\overline{\partial}(G_{\gra\grb})\\ 
   & \equiv & \ca d(G_{\gra\grb}^{-1}).\theta_{\gra}
                 + \overline{\partial}(G_{\gra\grb}),
\end{eqnarray*}
which proves (10) and yields, in turn, an $\ca$-linear connection
$\nabla$ on $\ce$. For the sake of completeness we outline the
construction of $\nabla$, referring for details to \cite{Mal3, V2}.
First, for each $\gra \in I$, we define the map
$\nabla^{\gra} : \ce|_{\ua}\lra \ce \otimes_{\ca}\Omega^1\,|_{\ua}$
by setting
\[
  \nabla^{\gra}(s) := \sum_{i=1}^{n}e^{\gra}_i \otimes
      \big(\partial(s_i^{\gra}) +
     \sum_{j=1}^n s_j^{\gra}\cdot \theta_{ij}^{\gra}\big),
\]
for every $s = \sum_{i=1}^n s_i^{\gra}\cdot e_i^{\gra} \in \ce(\ua)$
with $s_i^{\gra} \in \ca(\ua)$. Recall that $\theta_{\gra} \equiv
     (\theta_{ij}^{\gra})$, after the identifications (9).
The compatibility condition (10) implies that $\nabla^{\gra} = 
  \nabla^{\grb}$ on $\ca(\uab)$, hence we obtain a global
  connection $\nabla$.
\hfill $\Box$

\medskip
An immediate consequence of (11) is the following 

\begin{cor}    
If $\of: \cl \lra \ma$ is an isomorphism, then the connections
of $\cp$ are in bijective correspondence with the $\ca$-linear
connections of its associated vector sheaf $\ce$.
\end{cor}

\noindent
{\bf Example\ } Let $\ce$ be now a {\em given\/} vector sheaf of
rank $n$ with a local structure as in Paragraph 1.4. We denote
by $\cb$ the basis of topology on $X$ containing all the open
$V \sst X$ such that $V \sst \ua$, for some $\ua \in \cu$, and
consider the (complete) presheaf
$\cb \ni V \mapsto Iso_{\ca|_V}(\ca^n|_V, \ce|_V)$, where the
last space is the group of $\ca|_V$-linear isomorphisms. This
generates a principal sheaf $\cp(\ce) \equiv (\cp(\ce),\gla,X,
\pi)$, called {\em the sheaf of frames\/} of $\ce$.

We recall that (see \cite{V2}) there is a natural action of
$\gla$ on the right of $\cp(\ce)$ induced by the partial actions
\[
 Iso_{\ca|_V}(\ca^n|_V, \ce|_V)\times \gla (V) \lra
                      Iso_{\ca|_V}(\ca^n|_V, \ce|_V):
              (f,g) \mapsto f\cdot g \equiv f \circ g
\]
by employing, of course, the identification (2). The local
structure is described as follows: First we define the local
$\gla$-equivariant isomorphism
\[
  \Phi_{\gra}^V : \cp(\ce)(V) \lra \gla(V) :
            f \mapsto \psi_{\gra}^{-1} \circ f,                        
\]
for every open $V \in \ua$. Hence, varying $V$ in $\ua$, we obtain
an equivariant morphism
$\Phi_{\gra}: \cp(\ce)|_{\ua} \stackrel{\simeq}{\lra} \gla|_{\ua}$
and similarly for all $\gra \in I$.

The natural sections $\sigma_{\gra} \in \cp(\ce)(\ua)$, with respect
to ${\cal U}$, are now given by
\begin{equation}  
  \sigma_{\gra} := \Phi_{\gra}^{-1}(id|\ca^n(\ua)) = \psi_{\gra}.
\end{equation}

\medskip
The previous considerations lead now to 

\begin{cor}  
Every vector sheaf $\ce$ is associated with its principal sheaf of
frames $\cp(\ce)$, with respect to the trivial representation
of $\gla$. Hence, the $\ca$-linear connections
on $\ce$ correspond bijectively to the connections on $\cp(\ce)$.
\end{cor}

\noindent
{\bf Proof.} By the general construction discussed in Section
2, the vector sheaf, say, ${\cal F}$ associated with $\cp(\ce)$,
is generated by the presheaf
\[
  \cb \ni V \mapsto \cp(\ce)(V) \times \ca^n(V)/ \sim,
\]
defined by the trivial representation (: $\vf = id|\gla$). Though
we are restricted on a basis of topology, instead of the whole
topology of X, the final result remains unaffected. Following the
proof of \cite[Proposition 4.3]{V2}, for any $V \in \cb$ with
$V \sst \ua$, we consider the map
\[
   F_V : \cp(\ce)(V) \times \ca^n(V)/ \sim \;\;\;  \lra
              \;\; \ce(V) :[f,a] \mapsto f\circ a.
\]
We show that $F_V$ is a well defined bijection. Varying $V$ in
$\cb$, we obtain an isomorphism $F : {\cal F} \lra \ce$. It is
also an isomorphism of $\ca$-modules. Indeed, if we denote by
$\Psi_{\gra} : \ca^n|_{\ua} \lra {\cal F}|_{\ua} $ the isomorphisms
describing the local structure of ${\cal F}$, then (12) implies that
\[
  (F\circ \Psi_{\gra})(a) = F([\sigma_{\gra},a]) =
  \sigma_{\gra} \circ a = \psi_{\gra}(a);\quad a \in\ca^n(\ua),
\]
with $F$ and $\psi_{\gra}$ denoting now the induced morphisms between
sections. By the procedure used repeatedly so far, we see that
$F = \psi_{\gra} \circ \Psi_{\gra}^{-1}$. This, along
with the definition of the module operations on ${\cal F}$
(see Theorem 2), completes the claim about $F$. The rest of the
proof is clear. \hfill $\Box$

\begin{rem}
{\em In the previous Corollary we recover, by a different
approach, some of the results of \cite{V2}, notably Proposition 4.3
and Theorem 5.5.}
\end{rem}

\vspace{1.5cm}
\noindent
Efstathios Vassiliou\\
{\em Institute of mathematics\\
 University of Athens\\
 Panepistimiopolis\\
 Athens 157 84, Greece}

\medskip
\noindent
{\em E-mail address\/}:\verb= evassil@atlas.uoa.gr=
                        
\end{document}